\documentclass{conm-p-l}

      \title[Controlled Algebraic \textit{G}-theory, II]
            {Controlled Algebraic \textit{G}-theory, II}
     \author[Boris Goldfarb]{Boris Goldfarb}
    \address{Department of Mathematics and Statistics, SUNY, Albany, NY 12222}
      \email{bgoldfarb@albany.edu}
    \urladdr{http://www.albany.edu/~goldfarb/}
     \author[Timothy K.~Lance]{Timothy K.~Lance}
    \address{Department of Mathematics, Berkshire School, Sheffield, MA 01257}
      \email{tlance@berkshireschool.org}
       \date{\today}
     \thanks{The authors acknowledge support from the National Science Foundation.}

\usepackage{amsmath,amsthm,amssymb,amsxtra,amscd,enumerate,mathrsfs,graphicx,pifont}

\usepackage[all,cmtip]{xy}
  \SelectTips{cm}{10}\UseTips 

\AtBeginDocument{
\addtolength{\abovedisplayskip}{2pt}
\addtolength{\abovedisplayshortskip}{2pt}
\addtolength{\belowdisplayskip}{2pt}
\addtolength{\belowdisplayshortskip}{2pt}
}

\theoremstyle{plain}

\theoremstyle{definition}           

{\leqno{\normalshape ({\thedummy})} $$}
\numberwithin{equation}{dummy}

\theoremstyle{plain}

\theoremstyle{plain}
\newtheorem{Thm}{Theorem}[section]

\newtheorem{Cor}[Thm]{Corollary}
\newtheorem{Lem}[Thm]{Lemma}
\newtheorem{Prop}[Thm]{Proposition}

\theoremstyle{definition}
\newtheorem{Def}[Thm]{Definition}

\newtheorem{Ex}[Thm]{Example}
\newtheorem{Rem}[Thm]{Remark}

\theoremstyle{remark}
\newtheorem{Not}[Thm]{Notation}

\DeclareMathOperator{\A}{\mathbf{A}}

\DeclareMathOperator{\B}{\mathbf{B}}

\DeclareMathOperator{\C}{\mathbf{C}}

\DeclareMathOperator{\co}{\mathrm{co}}
\DeclareMathOperator{\coim}{coim}
\DeclareMathOperator{\coker}{coker}

\DeclareMathOperator{\D}{\mathbf{D}}

\DeclareMathOperator{\dist}{dist}
\DeclareMathOperator{\Dw}{\mathbf{D^{\boldsymbol{w}}}}

\DeclareMathOperator{\eU}{\mathbf{eU}}
\DeclareMathOperator{\E}{\mathbf{E}}

\DeclareMathOperator{\G}{\mathit{G}}

         \newcommand{\Gnc}{\G^{-\infty}}

\DeclareMathOperator{\im}{im}

\DeclareMathOperator{\K}{\mathit{K}}
         \newcommand{\Knc}{\K^{-\infty}}

\DeclareMathOperator{\Mod}{\mathbf{Mod}}

\DeclareMathOperator{\mU}{\mathbf{mU}}

         \newcommand{\subdot}{\boldsymbol{\cdot}}
\DeclareMathOperator{\supp}{supp}

\DeclareMathOperator{\U}{\mathbf{U}}

\DeclareMathOperator{\w}{\boldsymbol{w}}
\DeclareMathOperator{\bsv}{\boldsymbol{v}}

\DeclareMathOperator{\vD}{\mathbf{\boldsymbol{v}D}}
\DeclareMathOperator{\vDw}{\mathbf{\boldsymbol{v}D^{\boldsymbol{w}}}}
\DeclareMathOperator{\vE}{\mathbf{\boldsymbol{v}E}}

         \newcommand{\bfw}{\boldsymbol{w}}

\DeclareMathOperator{\vC}{\mathbf{\boldsymbol{v}C}}

\DeclareMathOperator{\wC}{\mathbf{\boldsymbol{w}C}}
\DeclareMathOperator{\wD}{\mathbf{\boldsymbol{w}D}}

\DeclareMathOperator{\Z}{\mathbf{Z}}

\newcommand{\pull}
{\!\!\! -\!\!\! -\!\!\! -\!\!\!}
\DeclareMathOperator*{\hocolimprep}{hocolim}
\newcommand{\hocolim}[1]%
{\hocolimprep_{\substack{- \pull \rightarrow \\ #1}} \, }
\DeclareMathOperator*{\colimprep}{colim}
\newcommand{\colim}[1]%
{\colimprep_{\substack{- \pull \rightarrow \\ #1}} \, }

\newcommand{\define}[1]{\textit{#1}}

\providecommand{\bysame}{\makebox[3em]{\hrulefill}\thinspace}

\newif\ifShowLabels
\ShowLabelstrue
\newcommand{\TeXref}[1]{
\marginpar{\scriptsize \texttt{#1}}}

\newcommand{\SecRef}[2]{\section{#1}\label{S:#2}%
\ifShowLabels \TeXref{{S:#2}} \fi}

\newcommand{\refT}[1]{\textup{\ref{T:#1}}}
\newcommand{\refL}[1]{\textup{\ref{L:#1}}}
\newcommand{\refD}[1]{\textup{\ref{D:#1}}}
\newcommand{\refC}[1]{\textup{\ref{C:#1}}}

\newcommand{\refP}[1]{\textup{\ref{P:#1}}}

\newenvironment{ThmRef}[1]%
{ \begin{Thm} \label{T:#1}
\ifShowLabels \TeXref{T:#1} \fi }%
{ \end{Thm} }
\newenvironment{DefRef}[1]%
{ \begin{Def} \label{D:#1}
\ifShowLabels \TeXref{D:#1} \fi }%
{ \end{Def} }
\newenvironment{LemRef}[1]%
{ \begin{Lem} \label{L:#1}
\ifShowLabels \TeXref{L:#1} \fi }%
{ \end{Lem} }
\newenvironment{CorRef}[1]%
{ \begin{Cor} \label{C:#1}
\ifShowLabels \TeXref{C:#1} \fi }%
{ \end{Cor} }
\newenvironment{RemRef}[1]%
{ \begin{Rem} \label{R:#1}
\ifShowLabels \TeXref{R:#1} \fi }%
{ \end{Rem} }
\newenvironment{PropRef}[1]%
{ \begin{Prop} \label{P:#1}
\ifShowLabels \TeXref{P:#1} \fi }%
{ \end{Prop} }
{ \begin{Ex} \label{E:#1}
\ifShowLabels \TeXref{E:#1} \fi }%
{ \end{Ex} }
\newenvironment{NotRef}[1]%
{ \begin{Not} \label{N:#1}
\ifShowLabels \TeXref{N:#1} \fi }%
{ \end{Not} }

\newenvironment{ThmRefName}[2]%
{ \begin{Thm} [#2]\label{T:#1}
\ifShowLabels \TeXref{T:#1} \fi }%
{ \end{Thm} }
\newenvironment{DefRefName}[2]%
{ \begin{Def} [#2]\label{D:#1}
\ifShowLabels \TeXref{D:#1} \fi }%
{ \end{Def} }
{ \begin{Lem} [#2]\label{L:#1}
\ifShowLabels \TeXref{L:#1} \fi }%
{ \end{Lem} }
{ \begin{Cor} [#2]\label{C:#1}
\ifShowLabels \TeXref{C:#1} \fi }%
{ \end{Cor} }
{ \begin{Rem} [#2]\label{R:#1}
\ifShowLabels \TeXref{R:#1} \fi }%
{ \end{Rem} }
{ \begin{Prop} [#2]\label{P:#1}
\ifShowLabels \TeXref{P:#1} \fi }%
{ \end{Prop} }
{ \begin{Ex} [#2]\label{E:#1}
\ifShowLabels \TeXref{E:#1} \fi  }%
{ \end{Ex} }

\ShowLabelsfalse

\begin{document}

\begin{abstract}
There are two established ways to introduce geometric control in the category of free modules---the bounded control and the continuous control at infinity.
Both types of control can be generalized to arbitrary modules over a noetherian ring and applied to study algebraic $K$-theory of infinite groups.
This was accomplished for bounded control in part I of the present paper and the subsequent work of G.~Carlsson and the first author, in the context of spaces of finite asymptotic dimension.
This part II develops the theory of filtered modules over a proper metric space with a good compactification.
It is applicable in particular to CAT(0) groups which do not necessarily have finite asymptotic dimension.
\end{abstract}

\maketitle

\tableofcontents

\SecRef{Introduction}{Intro}

Controlled algebra, and more specifically the algebraic $K$-theory of categories of controlled modules, is the foundation for recent computations of Loday's assembly map and the $K$-theory of new large classes of infinite groups \cite{aB:03,BFJR:04,gC:95,gCbG:04,gCeP:93,gCeP:97,dRrTgY:11}.
There are two established ways to introduce geometric control in the category of free modules.
The bounded control was introduced by Pedersen--Weibel \cite{ePcW:85,ePcW:89} and used by Gunnar Carlsson to study the
integral Loday assembly map \cite{gC:95}.  More recently, Carlsson and the first author applied similar ideas to
the question of surjectivity of the assembly maps for a large class of group rings \cite{gCbG:03,gCbG:02,gCbG:13}.
As part of that work, it was necessary to extend the constructions and the excision results to controlled modules that are not necessarily free.
This was done precisely under the name \textit{bounded $G$-theory} in \cite{gCbG:02}.

In this paper, we generalize the theory of continuously controlled geometric modules of Anderson--Connolly--Ferry--Pedersen \cite{dAfCsFeP:94}.
We do this with the view of extending the argument for the integral Borel Conjecture in $K$-theory to the settings which at this point can only be studied through continuous control.
The most important example of groups we have in mind are the CAT(0)-groups.

In sections 2 and 3, we will define exact categories of filtered $R$-modules over a metric space $M$
which are sensitive to global features of $M$.  The geometric control is imposed through the construction of a boundary $Y$ of $M$ so that the union of $M$ and $Y$ is a compact metric space $X$ which contains $M$ as an open dense subspace.  The space $X$ is usually referred to as the \textit{compactification of} $M$ by \textit{attaching the boundary} $Y$.

The nonconnective $K$-theory of our category associated to $X$, $Y$, and a given noetherian ring $R$ will be denoted $\Gnc (X,Y; R)$ or simply $\Gnc (X,Y)$ because the entire story is independent of a particular choice of $R$.
This theory is part of the commutative square
\[
\begin{CD}
\Knc (M,R) @>>> \Gnc (M,R) \\
@VVV @VVV \\
\Knc (X,Y;R) @>>> \Gnc (X,Y;R)
\end{CD}
\]
where $\Knc (M,R)$ is the nonconnective bounded $K$-theory of free $R$-modules of Pedersen--Weibel, $\Knc (X,Y;R)$ is the nonconnective delooping of the continuously controlled $K$-theory of Anderson et al., $\Gnc (M,R)$ is the bounded $G$-theory from \cite{gCbG:02}, and all maps of spectra are induced by exact inclusions between the underlying categories of modules.

The continuous control at infinity for free geometric modules, introduced by Anderson et al.~\cite{dAfCsFeP:94}, has been used to prove various versions of the integral Novikov conjecture \cite{gCeP:93,gCeP:97,bG:97,bG:98,bG:99,dR:04}.
The most important technical property of the continuously controlled $K$-theory of these geometric modules is the controlled excision, which allows to localize to subsets of the boundary and describes the continuously controlled $K$-theory in terms of these localizations.
The second goal of this paper is to prove an appropriate version of controlled excision for the continuously controlled $G$-theory.
This is accomplished in sections 4 and 5.
Suppose $\{ U, V \}$ is an open covering of $Y$.
If $\Gnc (X,Y; R)_{<U}$ stands for the $K$-theory of filtered modules over $M$ localized near $U$ in the appropriate sense that we specify later, the controlled excision is expressed by the homotopy pushout square of spectra
\[
\begin{CD}
\Gnc (X,Y; R)_{<U \cap V} @>>> \Gnc (X,Y; R)_{<U} \\
@VVV @VVV \\
\Gnc (X,Y; R)_{<V} @>>> \Gnc (X,Y; R)
\end{CD}
\]

To dispel the reader's suspicion that this paper is simply an exercise in combining the features of the two theories $\Gnc (M,R)$ and $\Knc (X,Y;R)$, we point out that the construction of $\Gnc (M,R)$ in \cite{gCbG:02} is far from a direct translation of Pedersen--Weibel \cite{ePcW:89}.
The challenge is in preserving the balance between the increasingly global nature of the objects and the need to localize to metric subspaces or subsets of the boundary in order to obtain excision theorems.
In many ways the outcome $\Gnc (X,Y; R)$ is more natural and the statements are simpler in this setting employing the continuous control.

The authors gratefully acknowledge support from the National
Science Foundation.

\SecRef{Continuous Control for Filtered Modules}{Sheaves}

Given a noetherian ring $R$ and a compactification of a metric space $M$, we will define a category of filtrations of left $R$-modules by subsets of $M$ and relate it to the category of geometric continuously controlled modules of Anderson et al.

\begin{DefRef}{MFIL}
Let $R$ be a left noetherian ring with unit and $\Mod(R)$ be the
category of left $R$-modules. For a general topological space $M$, the power set $\mathcal{P}(M)$ partially ordered by inclusion is a category
with subsets of $M$ as objects and unique morphisms $(S,T)$ when
$S \subset T$.

An
\define{$M$-filtered $R$-module} is a
functor
\[
F \colon \mathcal{P}(M) \longrightarrow \Mod (R)
\]
where
all structure maps $F(S,T)$ are
monomorphisms and $F(\emptyset) = 0$.
We will abuse notation by referring to $F(M)$ as $F$ when the filtration of $F(M)$ is clear from the context.

For any covariant functor $F \colon \mathcal{P}(M) \to \Mod (R)$ with $F(\emptyset) = 0$ there exists an associated
$M$-filtered $R$-module $F_M$ given by
$F_M (S) = \im F (S,M)$.

Given a filtered module $F$ and an arbitrary submodule $F'$,
the submodule has the \textit{standard $M$-filtration} given by $F' (S) = F(S) \cap F'$.
For example, for each subset $T \subset M$, the submodule $F(T)$ of $F$ has the canonical filtration and gives a filtered
module $F_T$ with $F_T (S) = F(S) \cap F(T)$.
\end{DefRef}

From now on we assume that $M$ is a proper metric space in the sense that all closed
metric balls in $M$ are compact.

Recall that a \textit{compactification} of $M$ is a compact metric space $X$ that contains the metric subspace $M$ as an open dense subset.
We also say that the compactification $X$ is obtained \textit{by attaching the boundary} $Y = X \setminus M$, or that $X$ is a compactification of $M$ \textit{by the boundary} $Y$.

\begin{DefRef}{ContCat}
Suppose $X$ is a compactification of $M$ by the boundary $Y$, and let $F$ and $G$ be two $M$-filtered modules.

An $R$-homomorphism $f \colon F \to G$
is \define{continuously controlled at} $y \in Y$
if for every neighborhood $U_y$ of $y$ in $X$
there is a neighborhood $V_y \subset U_y$ that satisfies
\begin{align}
f F(M \setminus U_y) &\subset G(M \setminus V_y), \notag \\
f F(V_y \cap M) &\subset G(U_y \cap M). \notag
\end{align}

Let $T$ be an open subset of $Y$.
An $R$-homomorphism $f \colon F \to G$
is \define{continuously controlled near} $T$ if for any open neighborhood $U$ of $T$ in $X$,
 there is another neighborhood $V$  of $T$ that satisfies
\begin{align}
f F(M \setminus U) &\subset G(M \setminus V), \tag{1} \\
f F(V \cap M) &\subset G(U \cap M). \tag{2}
\end{align}

If $f$ is continuously controlled near all open subsets $T \subset Y$, we call $f$ \define{continuously controlled}.
\end{DefRef}

\begin{PropRef}{IOUY}
If an $R$-homomorphism $f$ between $M$-filtered modules is continuously controlled then it is continuously controlled at all points of $Y$.
\end{PropRef}

\begin{proof}
Consider $y$ in $Y$ and an open neighborhood $U_y$.  If $T = U_y \cap Y$ then there is $V$ with $V \cap Y = T$ and the required properties.  We can choose $V_y = V$.
\end{proof}

The converse to Proposition \refP{IOUY} is not true without additional assumptions on the compactification $(X,Y)$ and the filtrations of $F$ and $G$.

\begin{DefRef}{ContCat2}
Suppose an $R$-homomorphism $f \colon F \to G$
is continuously controlled near $T$, where $T$ is an open subset of $Y$.
That is, for every $U$ there is $V$ as in Definition \refD{ContCat} so that conditions (1) and (2) hold.

If, in addition,
\begin{align}
G(M \setminus U) \cap f F(M) &\subset f F(M \setminus V), \tag{3} \\
G(V \cap M) \cap f F(M) &\subset f F(U \cap M), \tag{4}
\end{align}
we will say that $f$ is \define{continuously bicontrolled near} $T$.
If $f$ is continuously bicontrolled near all open subsets $T \subset Y$, we call $f$ \define{continuously bicontrolled}.
\end{DefRef}

\begin{NotRef}{LE}
If all four of the conditions (1) through (4) are satisfied, we
denote this relationship by $V \le_{f} U$ when $T \subset Y$ is
understood.
In the case when $V \le_{f} V$ for all $V$, we say that $f$ is \textit{static at infinity}.
\end{NotRef}

We will briefly recall the construction of the geometric $R$-modules over $(X,Y)$ and explain its relation to the filtered $R$-modules.
The construction in \cite{dAfCsFeP:94} is given for arbitrary additive categories but we specialize to the case of finitely generated $R$-modules.

\begin{DefRefName}{KOIUN}{Anderson--Connolly--Ferry--Pedersen \cite{dAfCsFeP:94}}
Let $(X,Y)$ be a compactification of $M$ as before, and let $R$ be a ring.

The category $\mathcal{C} (X,Y)$ is defined as follows.
An object of $\mathcal{C} (X,Y)$ is a collection $A$ of  finitely generated free $R$-modules
$\{ A_x \mid x \in X \setminus Y \}$
with the property that for each compact subset $K$ of $M = X \setminus Y$ the set $\{ x \in K \mid A_x \ne 0 \}$ is finite.
A morphism $f \colon A \to B$ in $\mathcal{C} (X,Y)$ is a collection of $R$-module homomorphisms $\{ f^x_y \colon A_x \to B_y \}$ such that (a) for all $x$ the set $\{ y \mid f^x_y \ne 0 \}$ is finite and (b) for every point $z \in Y$ and every neighborhood $U$ of $z$ in $X$ there is a neighborhood $V$ of $z$ such that $f^x_y = 0$ whenever $x$ in $V$ and $y$ is not in $U$.
Such morphisms are called \textit{continuously controlled at infinity}.

The additive structure of $\mathcal{C} (X,Y)$ is given by defining $A \oplus B$ as the collection $\{ A_x \oplus B_x \}$.
\end{DefRefName}

To see that objects of $\mathcal{C} (X,Y)$ correspond to certain natural filtered modules, given $A$ as above we define the $R$-module $F_A = \bigoplus_x A_x$ and its $M$-filtration by $F_A (S) = \bigoplus_{x \in S} A_x$.
The structure maps are the inclusions of direct summands.
It is easy to see that $F_A$ is has the property that $F_A (K)$ is a finitely generated $R$-module for a compact subset $K$.  All morphisms $f \colon A \to B$ which are continuously controlled at infinity are in fact continuously controlled morphisms $f \colon F_A \to F_B$.
This means that $\mathcal{C} (X,Y)$ is a full subcategory of $M$-filtered $R$-modules and continuously controlled homomorphisms.

Again, not every morphism from $\mathcal{C} (X,Y)$ is going to be continuously bicontrolled.
However, inclusions of direct summands are bicontrolled.
The importance of that more restrictive condition is evident in the introduction of a new exact structure on filtered modules, which is the next goal.
We will see that $\mathcal{C} (X,Y)$ is in fact an exact subcategory of the new exact category.

\begin{LemRef}{CatOK}
Let $f_1 \colon F \to G$, $f_2 \colon G \to H$ be continuously controlled homomorphisms between $M$-filtered modules with respect to some compactification $X$ of $M$.
and $f_3 = f_2 \circ f_1$.
\begin{enumerate}
\item If $f_1$, $f_2$ are continuously bicontrolled morphisms
and either $f_1 \colon F(M) \to G(M)$ is an epi or $f_2 \colon
G(M) \to H(M)$ is a monic, then $f_3$ is also continuously
bicontrolled.
\item If $f_1$, $f_3$ are continuously bicontrolled
and $f_1$ is epic then $f_2$ is also continuously bicontrolled; if
$f_3$ is only continuously controlled then $f_2$ is also continuously
controlled.
\item If $f_2$, $f_3$ are continuously bicontrolled
and $f_2$ is monic then $f_1$ is also continuously bicontrolled; if
$f_3$ is only continuously controlled then $f_1$ is also continuously
controlled.
\end{enumerate}
\end{LemRef}

\begin{proof}
To see, for example, part (2) of the statement, let $f_1$,
$f_3$ be continuously bicontrolled homomorphisms
so that $f_1$ is surjective,
and let $U \subset X$ be any neighborhood of $y \in Y$.
If $V$, $V'$ are neighborhoods of $y$ with
$V' \le_{f_1} U$ and $V \le_{f_3} V'$ then $V \le_{f_2} U$
as required.
\end{proof}

\begin{DefRef}{HUTWS}
A \define{pseudoabelian} category is an
additive category with kernels and cokernels and the property that
for every morphism $f$, the canonical map $\pi \colon \coim
(f) \to \im (f)$ is both monic and epic. In particular,
$f$ factors canonically as $f = m  e$ where $m = \im (f)
= \ker (\coker f)$ is monic and $e$, the composition of the
epic $\coim (f) = \coker (\ker f)$ and $\pi$, is epic. The
category is
\define{abelian} if, in addition, the map $\pi$ is an isomorphism.
\end{DefRef}

The most common way pseudoabelian features come up in a category
is via the property that every monic is a kernel which a pseudoabelian category may or
may not have. This is generally not true in categories of interest in this paper.

Let $R$ be a fixed noetherian ring and $X$ be a compactification of a proper metric space $M$ by the boundary $Y$.

\begin{DefRef}{PlainU}
We will call an $M$-filtered
$R$-module $F$ \textit{locally finite} if
the subcategory $\mathcal{B}(M)$ of
bounded subsets maps to the subcategory of finitely generated submodules of $F(M)$.

The objects of the category $\U (X,Y)$ are the locally finite $M$-filtered $R$-modules.
The morphisms in  $\U (X,Y)$ are the continuously controlled homomorphisms.
\end{DefRef}

\begin{ThmRef}{ContAb}
The category $\U (X,Y)$ is a pseudoabelian category.
\end{ThmRef}

\begin{proof}
The additive properties are inherited from $\Mod (R)$. In particular, the
biproduct is given by the filtration-wise operation
$(F \oplus G)(S) = F(S) \oplus G(S)$.
Notice also that $\pi$ is monic or epic in
$\U (X,Y)$ if and only if it is such in $\Mod (R)$.

The kernel $K$ of an arbitrary continuously controlled
morphism $f \colon F \to G$ in $\Mod (R)$ has the
standard filtration
$K(S) = K \cap F(S)$
which gives the kernel of $f$.
Since $R$ is noetherian and $F$ is locally finite, $K$ is locally finite.
The canonical monic
$\kappa \colon K \to F$ is stable at infinity and therefore continuously bicontrolled.
It follows from
part (3) of Lemma \refL{CatOK} that $K$ has the universal properties
of the kernel in $\U (X,Y)$.

Similarly, let $I$ be the
$M$-filtration of the image of $f$ in $\Mod (R)$ by $I(S) = \im (f) \cap G(S)$.
Then there is a filtered module $C$ given by
$C(S) = G(S)/I(S)$
for $S \subset M$.
The module $C(M)$ is the cokernel of $f$ in
$\Mod (R)$.
As in Definition \refD{MFIL}, there is an $M$-filtered module $C_M$ associated
to $C$ given by $C_M (S) = \im C (S,M)$.
Since $R$ is noetherian and $G$ is locally finite, both $I$ and $C$ are locally finite.
The canonical homomorphism
$\sigma \colon G(M) \to C(M)$ gives a continuously bicontrolled morphism
$\sigma \colon G \to C_M$ which is stable at infinity since
$\im (\sigma  G(S,M)) = C_M (S)$.
This, in
conjunction with part (2) of Lemma \refL{CatOK}, also verifies the
universal properties of $C_M$ and $\sigma$ in $\U (X,Y)$.
\end{proof}

\begin{LemRef}{FutRef}
An isomorphism in $\U (X,Y)$ is continuously bicontrolled.
A morphism $f$ in $\U (X,Y)$ is continuously bicontrolled if and only
if it is balanced, that is, the canonical map $\pi \colon \coim (f) \to \im (f)$
is an isomorphism in $\U (X,Y)$.
\end{LemRef}

\begin{proof}
If $f^{-1}$ is the inverse to $f$ and $V \le_f U$ then $U \le_{f^{-1}} V$.
For instance, since $f$ is epic, condition (3) from Definition \refD{ContCat} is equivalent to
\[
G(M \setminus U) \subset f F(M \setminus V)
\]
which is equivalent to
\[
f^{-1} G(M \setminus U) \subset F(M \setminus V),
\]
which is condition (1) for $f^{-1}$.
The other six conditions for $f$ and $f^{-1}$ are similarly pairwise equivalent.
\end{proof}

\begin{ThmRef}{UXYex}
The continuously bicontrolled injections $\mU (X,Y)$ and
the continuously bicontrolled surjections $\eU (X,Y)$ are the
admissible monomorphisms and admissible epimorphisms for an exact structure in $\U (X,Y)$.
\end{ThmRef}

\begin{proof}
Formally the constructions in the proof are identical to those in the proof of Theorem 2.13 in \cite{gCbG:02}.
The verification of their properties is quite different.
We need to verify Quillen's axioms for the family $\mathcal{E}$
of exact sequences
\[
E^{\subdot} \colon \quad E' \xrightarrow{\ i \ } E \xrightarrow{\ j \ } E''
\]
where $i$ is an admissible monomorphism and $j$ is an admissible epimorphism.

(1) By Lemma \refL{FutRef}, any short exact sequence $F^{\subdot}$ isomorphic to
some $E^{\subdot} \in \mathcal{E}$ is also in $\mathcal{E}$. A
split exact sequence $E' \to E' \oplus E'' \to E''$ is in
$\mathcal{E}$.  Also $i = \ker (j)$ and $j = \coker (i)$ in any exact sequence
$E^{\subdot} \in \mathcal{E}$.

(2) The collections $\mU (X,Y)$ and $\eU (X,Y)$ are
closed under composition by part (1) of Lemma \refL{FutRef}.
It remains to check that $\mU (X,Y)$ and $\eU (X,Y)$ are
closed under base and cobase changes.

Given an exact sequence $E^{\subdot} \in \mathcal{E}$
and a morphism $f \colon A \to E''$, there is a base
change diagram
\[
\begin{CD}
E' @>>> E \times_{f} A @>{j'}>> A \\
@V{=}VV @VV{f'}V @VV{f}V \\
E' @>>> E @>{j}>> E''
\end{CD}
\]
where the module
$E \times_{f} A$
is the fiber product of $j$ and $f$.
Since $R$ is noetherian and both $E$ and $A$ are locally finite, the
standard filtration given by
\begin{align}
(E \times_{f} A)(S) & = \{ (e,a) \in E \times A \mid j(e) =
f(a), e\in E(S), a \in A(S) \} \notag \\
& = E
\times_{f} A \cap \big( E(S) \times A(S) \big) \notag
\end{align}
is locally finite.
To check that the induced
epimorphism $j'$ is continuously controlled, given a neighborhood $U$ of $y \in Y$ we need a neighborhood $V$ such that
\[
j' (E \times_{f} A)(M \backslash U) \subset A (M \backslash V).
\]
If $x = (e,a) \in (E \times_{f} A)(M \backslash U)$ then by definition $x \in A (M \backslash U) \times E (M \backslash U)$, so $j' (x) = a \in A (M \backslash U)$.  So we can choose $V = U$.  The same choice verifies the other control property.  This shows that $j'$ is stable at infinity.
To check that $j'$ is continuously bicontrolled, use the fact that $f$ is continuously controlled to choose $U''$ for the given $U$ so that
$f A (M \backslash U) \subset E'' (M \backslash U'')$.
Using that $j$ is continuously bicontrolled, choose $V$ such that
\[
E'' (M \backslash U'') \cap j(E) \subset jE (M \backslash V).
\]
Given
$a \in A (M \backslash U) \cap j' (E \times_{f} A)$,
we have $f(a) \in E'' (M \backslash U'')$.
Since $f(a) \in j(E)$, there is $\overline{e} \in E (M \backslash V)$ with $j(\overline{e}) = f(a)$.
Now
$(\overline{e},a) \in (E \times_{f} A)(M \backslash V)$
and $j' ((\overline{e},a)) = a$.  The other bicontrol property follows similarly.
It also follows easily that $i'$ is continuosly bicontrolled.
So the class
of admissible epimorphisms is closed under base change by
arbitrary morphisms in $\U (X,Y)$.

The pushout in the cobase
change diagram
\[
\begin{CD}
E' @>{i}>> E @>{j}>> E'' \\
@V{g}VV @VVV @VV{=}V \\
B @>{i'}>> E \oplus_{g} B @>{j'}>> E''
\end{CD}
\]
is the cokernel of the continuously controlled inclusion $E' \to E
\oplus B$ sending $e'$ to $(-i(e'), g (e'))$.
The induced injection $i' \colon B
\to E \oplus_{g} B$ is clearly continuously bicontrolled and $j'$ is static at infinity.

(3) It is known that the third exactness axiom of Quillen follows from the first two, cf.~\cite{bK:96}.  This completes the proof.
\end{proof}

\begin{PropRef}{SSSST}
The exact structure $\mathcal{E}$ in $\U (X,Y)$ consists of sequences isomorphic to those
\[
E^{\subdot} \colon \quad E' \xrightarrow{\ i \ } E \xrightarrow{\ j \ } E''
\]
which possess restrictions
\[
E^{\subdot} (S) \colon \quad E' (S) \xrightarrow{\ i \ } E (S) \xrightarrow{\ j \ } E'' (S)
\]
for all subsets $S \subset M$, and each $E^{\subdot} (S)$ is an exact sequence in $\Mod (R)$.
\end{PropRef}

\begin{proof}
Compare to Proposition 2.14 in \cite{gCbG:02}.
\end{proof}

\SecRef{Continuously Controlled $G$-theory}{CCGT}

\begin{DefRef}{KLPOI}
Given a subset $Z$ of $Y$, a submodule $H$ of a filtered module $F$ is \textit{supported near $Z$}, written $\supp (H) \subset Z$,  if there is a closed subset $C$ of $X$ with $H \subset F(C \cap M)$ and $C \cap Y \subset Z$.

A filtered module $F$ is \textit{insular at infinity} if, given any pair of subsets $T_1$ and $T_2$ of $Y$ and a submodule $H$ of $F$,
whenever
\[
\supp H \subset T_1, \ \supp H \subset T_2
\]
then
\[
\supp H \subset T_1 \cap T_2.
\]

A filtered module $F$ is \textit{split at infinity} if, given any pair of subsets $U_1$ and $U_2$ of $Y$ and a submodule $H$ of $F$ with
\[
\supp (H) \subset U_1 \cup U_2,
\]
there are two submodules $H_1$ and $H_2$ of $H$ such that
$H \subset H_1 + H_2$
and
\[
\supp H_1 \subset U_1, \ \supp H_2 \subset U_2.
\]
\end{DefRef}

One easy observation is that checking splitting at infinity reduces to checking it for closed subsets of $Y$, because of the definition of support.

\begin{LemRef}{GGGGG}
Let $T$ be a subset of $Y$.

$\mathrm{(a)}$  If $f \colon F \to G$ is continuously controlled at infinity and $\supp (H) \subset T$ for some submodule $H \in F$ then $\supp f(H) \subset T$.

$\mathrm{(b)}$  If $f$ is continuously bicontrolled at infinity and $\supp (H) \subset T$ for some submodule $H \subset \im (f)$ then there is a submodule $P \in F$ such that $\supp (P) \subset T$ and $H \subset f(P)$.

$\mathrm{(c)}$  If $\supp (H_1) \subset T$ and $\supp (H_2) \subset T$ then $\supp (H_1 + H_2) \subset T$.

$\mathrm{(d)}$  If $H' \subset H$ and $\supp (H) \subset T$ then $\supp (H') \subset T$.
\end{LemRef}

\begin{proof}
(a) Suppose $\supp (H) \subset T$, so there is a closed $C$ with $H \subset F(C \cap M)$ and $C \cap Y \subset T$.
Now $U = X \backslash C$ is an open neighborhood of the open set $T' = Y \backslash C$, and there is another open neighborhood $V \subset X \backslash C$ of $T'$ such that $V \cap Y = T'$ and $f(H) \subset G(M \backslash V)$.
Since $(X \backslash V) \cap Y = C \cap Y$, we have $\supp f (H) \subset T$.
The statement (b) is proved similarly.

(c) Suppose $C_1$ and $C_2$ are closed subsets of $X$ such that $H_1 \subset F(C_1 \cap M)$ and $H_2 \subset F(C_2 \cap M)$, and $C_i \cap Y \subset T$ for $i = 1$, $2$.  Then
$H_1 + H_2 \subset F((C_1 \cup C_2) \cap M)$,
and $(C_1 \cup C_2) \cap Y \subset T$.

(d) is evident.
\end{proof}

\begin{ThmRef}{lninpres}
$\mathrm{(1)}$  Insular at infinity objects are closed under exact extensions in
$\U (X,Y)$, that is, if
\[
E' \longrightarrow E \longrightarrow E''
\]
is an exact sequence in $\U (X,Y)$, and $E'$, $E''$ are insular at infinity, then
$E$ is insular at infinity.

$\mathrm{(2)}$  If $E$ in the exact extension above is insular at infinity then $E'$ is insular at infinity.

$\mathrm{(3)}$  Split at infinity objects are closed under exact extensions.

$\mathrm{(4)}$  If $E$ in the exact extension above is split at infinity then $E''$ is split at infinity.

$\mathrm{(5)}$  If $E$ is both insular at infinity and split at infinity then $E''$ is insular at infinity if and only if $E'$ is split at infinity.
\end{ThmRef}

\begin{proof}
Let
\[
E' \xrightarrow{\ f \ } E \xrightarrow{\ g \ } E'' \tag{$\ast$}
\]
be an exact sequence in $\U (X,Y)$.

(1) Since $E''$ is insular at infinity, for any pair
of subsets $T_1$ and $T_2$ of $Y$ and any submodule $H$ of $E$ with $\supp (H) \subset T_1$ and $\supp (H) \subset T_2$, we have
$\supp g (H) \subset T_1 \cap T_2$
using part (a) of Lemma \refL{GGGGG}.

Using part (b) of Lemma \refL{GGGGG}, choose a submodule $P$ of $E$ such that $g (H) \subset g(P)$ and $\supp (P) \subset T_1 \cap T_2$.
By part (c) of Lemma \refL{GGGGG}, $\supp (H + P) \subset T_1$ and $\supp (H + P) \subset T_2$, so $\supp \big( fE' \cap (H + P) \big) \subset T_1$ and $\supp \big( fE' \cap (H + P) \big) \subset T_2$.
Since $f$ is injective, part (b) shows that $\supp f^{-1} \big( fE' \cap (H + P) \big) \subset T_1$ and $\supp f^{-1} \big( fE' \cap (H + P) \big) \subset T_2$, so
\[
\supp f^{-1} \big( fE' \cap (H + P) \big) \subset T_1 \cap T_2
\]
as $E'$ is insular at infinity.  By part (a),
\[
\supp \big( fE' \cap (H + P) \big) \subset T_1 \cap T_2.
\]
Now $H \subset \big( fE' \cap (H + P) \big) + P$, so $\supp (H) \subset T_1 \cap T_2$ by part (c) of Lemma \refL{GGGGG}.

(2) Suppose $H' \subset E'$ and $\supp (H') \subset T_1$ and $\supp (H') \subset T_2$.
Using part (a) of Lemma \refL{GGGGG} and the fact that $E$ is insular at infinity, we have $\supp f(H') \subset T_1 \cap T_2$.
By part (b), $\supp (H') \subset T_1 \cap T_2$.

(3) Given a submodule $H \subset E$ with $\supp (H) \subset U_1 \cup U_2$, by part (a) of Lemma \refL{GGGGG} we have the submodule $gH \subset E''$ with $\supp g(H) \subset U_1 \cup U_2$.  $E''$ is split at infinity, so there are submodules $H_1^g$ and $H_2^g$ of $E''$ such that $\supp H_i^g \subset U_i$ for $i=1$, $2$, and $H \subset H_1^g + H_2^g$.

Since $g$ is surjective, we can find submodules $P_i$ of $E$ such that $\supp (P_i) \subset U_i$ and $H_i^g \subset g(P_i)$.
If $x \in H$ and we write $g(x) = y_1 + y_2$ with $y_i \in H_i^g$, this provides $x_i \in P_i$ such that $g(x_i) = y_i$.
Then $g(x- x_1 - x_2) = 0$, so $H \subset P_1 + P_2 + \ker (g)$.
By exactness, $\ker (g) = \im (f)$, so by part (c) we have $\supp (H + P_1 + P_2) \subset U_1 \cup U_2$.
By part (d), the submodule
\[
K = (H + P_1 + P_2) \cap \ker (g)
\]
has $\supp (K) \subset U_1 \cup U_2$.
Then we have $f^{-1} (K) \subset Q_1 + Q_2$ in $E'$ with $\supp (Q_i) \subset U_i$ and so $\supp (fQ_i) \subset U_i$.
Now
\[
H \subset P_1 + P_2 +f(Q_1 + Q_2) = (P_1 + fQ_1) + (P_2 + fQ_2),
\]
and $\supp (P_i + fQ_i ) \subset U_i$.

(4) Consider $H'' \subset E''$ with $\supp (H'') \subset U_1 \cup U_2$.  Then there is $P \subset E$ such that $H'' \subset f(P)$ and $\supp (P) \subset U_1 \cup U_2$.  Let $P = P_1 + P_2$ be the splitting decomposition, then $H'' \subset f(P_1) + f(P_2)$ and $\supp f(P_i) \subset U_i$.

(5) We will show that if $E'$ is split at infinity and $E$ is insular at infinity in the exact sequence ($\ast$), then $E''$ is insular at infinity.
The proof of the fact that $E''$ is insular at infinity and $E$ split at infinity implies $E'$ split at infinity is left to the reader.

Given $H \subset E''$ with $\supp (H) \subset T_1$ and $\supp (H) \subset T_2$, there are submodules $P_1$ and $P_2$ of $E$ such that $H_i \subset g(P_i)$ and $\supp (P_i) \subset T_i$ for $i=1$, $2$.

Suppose $y \in H$.
For $x_1 \in P_1$ and $x_2 \in P_2$ such that $g(x_1) = g(x_2) =y$, $k = x_1 - x_2$ is an element of $(P_1 + P_2) \cap \ker (g) = (P_1 + P_2) \cap \im (f)$.
Notice that $\supp (P_1 + P_2) \subset T_1 \cup T_2$, so there is a submodule $Q$ of $E'$ with $\supp (Q) \subset T_1 \cup T_2$
and $(P_1 + P_2) \cap \im (f) \subset f(Q)$.
Find $z \in Q$ such that $f(z) = k$.
Since $E'$ is split at infinity, there are submodules $Q_1$ and $Q_2$ of $E'$ with $\supp (Q_i) \subset T_i$ such that $Q \subset Q_1 + Q_2$.
So $z = z_1 + z_2$ for some $z_i \in Q_i$.
Now $k = f(z_1) + f(z_2)$ where $f(z_i) \in f(Q_i)$ and $\supp f(Q_i) \subset T_i$.
Therefore, $x_1 = x_2 + k = x_2 + f(z_1) + f(z_2)$.
So
\[
x = x_1 - f(z_1) = x_2 + f(z_2)
\]
is an element of both $P_i + f(Q_i)$, $i=1$, $2$, and $\supp (P_i + fQ_i) \subset T_i$.
We found $x$ in
\[
P = (P_1 + f(Q_1)) \cap (P_2 + f(Q_2))
\]
such that $g(x) = y$.
Notice that $\supp (P) \subset T_1 \cap T_2$ because $E$ is insular at infinity.
So $H \subset g(P)$ and $\supp (gP) \subset T_1 \cap T_2$.
By part (d) of Lemma \refL{GGGGG}, $\supp (H) \subset T_1 \cap T_2$.
\end{proof}

\begin{DefRef}{OPIML}
The category $\C' (X, Y)$ is the full subcategory of $\U (X, Y)$ on the filtered modules which are split and insular at infinity.
\end{DefRef}

\begin{ThmRef}{KAIR}
$\C' (X, Y)$ is closed under exact extensions in $\U (X, Y)$.
It is, therefore, an exact category.
\end{ThmRef}

\begin{proof}
It is clear that $\C' (X, Y)$ is closed under isomorphisms in $\U (X, Y)$.
Closure under exact extensions follows from parts (1) and (3) of Theorem \refT{lninpres}.
\end{proof}

Let $\mathcal{C} (X,Y)$ be the continuously controlled category of geometric $R$-modules of Anderson--Connolly--Ferry--Pedersen \cite{dAfCsFeP:94}.
It is an additive category and therefore has the split exact structure.

\begin{PropRef}{THYU}
$\mathcal{C} (X,Y)$  is an exact subcategory of $\C' (X,Y)$.
\end{PropRef}

\begin{proof}
The objects of $\mathcal{C} (X,Y)$ can be described as
the $M$-filtered modules with
is
\[
F(S) = \bigoplus_{x \in S} F(x),
\]
and the
structure maps that are the inclusions of direct summands.
It remains to check that the converse of Proposition \refP{IOUY} is true for morphisms between objects from $\mathcal{C} (X,Y)$,
that is, the homomorphisms controlled at all $y$ in $Y$ are controlled near all open subsets of $Y$ if both domain and target are in $\mathcal{C} (X,Y)$.  This is left to the reader.
\end{proof}

\begin{RemRef}{THYUPO}
Another consequence of Theorem \refT{lninpres} is a procedure for constructing interesting examples of filtered modules which are not geometric modules.
These can be obtained as kernels or cokernels of continuously bicontrolled homomorphisms.
A class of such homomorphisms is available as idempotents of geometric modules.
Indeed, under idempotents elements of modules can be pulled back identically.
For example, in the image of an idempotent which is not an inclusion onto a direct summand, the filtration
includes non-free projective $R$-modules.
\end{RemRef}

Next we would like to recall the definition of the boundedly controlled category $\B (M)$ from \cite{gCbG:02}.

\begin{DefRef}{RealBCprelim}
Given a proper metric space $M$, let $F$ be a locally finite $M$-filtered module, cf.~Definition \refD{PlainU}.

A filtered module $F$ is \textit{lean} or $D$-\textit{lean} if there is a number
$D \ge 0$ such that for every subset $S$ of $M$
\[
F(S) \subset \sum_{x \in S} F(B_D (x)),
\]
where $B_D (x)$ is the metric ball of radius $D$ centered at $x$.

A filtered module $F$ is \textit{insular} or $d$-\textit{insular} if there is a
number $d \ge 0$ such that
\[
F(T) \cap F(U) \subset F(T[d] \cap U[d])
\]
for every pair of subsets $T$, $U$ of $M$.

The \textit{boundedly controlled morphisms} between $M$-filtered modules are
$R$-linear homomorphisms $\phi \colon F_1 \to F_2$ for which there exists a number $b \ge 0$
such that the image
$\phi (F_1 (S))$ is contained in the submodule $F_2 (S [b])$ for
all subsets $S \subset M$. Here $S[b]$ stands for the
metric $b$-enlargement of $S$ in $M$.

Let $\B' (M)$ be the category of locally finite lean insular objects and bounded morphisms.

We say that a boundedly controlled morphism $\phi$ as above is in fact \textit{boundedly bicontrolled} if in addition to the inclusions
$\phi (F_1 (S)) \subset F_2 (S [b])$ one also has the inclusions $\im (\phi) \cap F_2 (S) \subset \phi (F_1 (S[b]))$ for all subsets $S$ of $M$.
\end{DefRef}

It was shown in \cite{gCbG:02} that $\B' (M)$ is an exact category with respect to the collection of exact sequences where the admissible morphisms are the boundedly bicontrolled injections and surjections.

\begin{DefRef}{SmallInf}
Given a pair $(X,Y)$ as before, the metric on $M = X \setminus Y$
is called \define{small at infinity} if for each $y \in Y$, each $D
\ge 0$, and a neighborhood $U \subset X$ of $y$, there is a
neighborhood $V \ni y$ such that for $x \in V$ one has $x[D] \subset U$.
In this case, we will also say that the compactification $(X,Y)$ is
\define{small at infinity}.
\end{DefRef}

\begin{PropRef}{THYU2}
Suppose $(X,Y)$ is small at infinity with respect to the metric on $M = X \setminus Y$.
Then the category $\B' (M)$ is an exact subcategory of $\C' (X,Y)$.
\end{PropRef}

\begin{proof}
It is easy to see that lean implies split at infinity and insular implies insular at infinity.
A boundedly controlled homomorphism is continuously controlled at infinity because the metric is small at infinity.
\end{proof}

\begin{CorRef}{KLSDA}
If the compactification $(X,Y)$ is small at infinity, there is a commutative square of exact inclusions
\[
\begin{CD}
\mathcal{B} (M) @>{\iota}>> \B' (M) \\
@V{\kappa}VV @VV{\overline{\kappa}}V \\
\mathcal{C} (X,Y) @>{\overline{\iota}}>> \C' (X,Y)
\end{CD}
\]
\end{CorRef}

\begin{proof}
Exactness of $\iota$ is known from \cite{gCbG:02}.
Exactness of $\kappa$ is known from \cite{dAfCsFeP:94}.
Exactness of $\overline{\kappa}$ follows from the smallness at infinity condition as in Proposition \refP{THYU2}.
For exactness of $\overline{\iota}$ one needs Proposition \refP{IOUY} to see that the morphisms in $\mathcal{C} (X,Y)$ are
continuously controlled in $\C (X,Y)$.
\end{proof}

\begin{DefRefName}{strict}{Definition 2.19 of \cite{gCbG:02}}
An object $F$ of $\B' (M)$ is called $\ell$-\textit{strict} or simply \textit{strict} if there exists an order preserving function
$\ell \colon \mathcal{P}(M) \to [0, \infty)$
such that for every subset $S$ of $M$ the subobject
$F_{S} = F(S)$
is
locally finite, $\ell_{S}$-insular and $\ell_{S}$-lean with respect to the
standard filtration
$F_{S} (U) = F_{S} \cap F(U)$.

The bounded category $\B (M)$ is defined in \cite{gCbG:02} as the full exact subcategory of $\B' (M)$ on filtered modules isomorphic to strict objects.
\end{DefRefName}

The following is a relaxation of the condition in the last definition.

\begin{DefRef}{graded}
An $M$-filtered module $F$ is called $\ell$-\textit{graded} or simply \textit{graded} if there exists an order preserving function
\[
\ell \colon \mathcal{P}(X) \longrightarrow [0, \infty)
\]
such that for every subset $S$ of $M$ there is a subobject $\widehat{F}_S$ such that
$\widehat{F}_S$
is
locally finite, $\ell_{S}$-insular and $\ell_{S}$-lean, and
\[
F(S) \subset \widehat{F}_S \subset F(S[\ell_S]).
\]
\end{DefRef}

Clearly, a strict object is graded.
It is easy to see that an object isomorphic to a strict object is also graded.
The graded objects are closed under isomorphisms and form an exact category $\B'' (M)$ containing $\B (M)$ as an exact subcategory.
In fact, the theory in \cite{gCbG:02} has an analogue entirely in terms of graded $M$-modules.
We don't require it in this paper.
We simply point out that the next definition is best viewed as an analogue of the graded filtered module rather than a strict module or isomorphic to strict.

\begin{DefRef}{graded2}
An object $F$ of $\C' (X,Y)$ is called $(X,Y)$-\textit{graded} or simply \textit{graded} if for each open subset $V \subset X$ there are
\begin{enumerate}
\item a submodule $F^V$ of $F$ which is locally finite, split and insular at infinity and a subset $V' \supset V$ such that $V' \cap Y = V \cap Y$ which satisfy
\[
F(V \cap M) \subset F^V \subset F(V' \cap M),
\]
\item a submodule $F_V$ of $F$ which is locally finite, split and insular at infinity and another subset $V'' \subset V$ such that $V'' \cap Y = V \cap Y$ which satisfy
\[
F(V'' \cap M) \subset F_V \subset F(V \cap M).
\]
\end{enumerate}

We define the continuously controlled category $\C (X,Y)$ to be the full subcategory of $\C' (X,Y)$ on the graded objects.
\end{DefRef}

\begin{PropRef}{POIYVB}
The category $\C (X,Y)$ is an exact category.  The inclusion of $\C (X,Y)$ in $\C' (X,Y)$ is exact.
\end{PropRef}

\begin{proof}
Given an exact sequence in $\C' (X,Y)$
\[
F \xrightarrow{\ f \ } G \xrightarrow{\ g \ } H,
\]
we assume that
$F$ and $H$ are graded modules in $\C (X,Y)$ with the associated gradings ${F}^{\bullet}$, ${F}_{\bullet}$, ${H}^{\bullet}$, and ${H}_{\bullet}$.

To define a grading for $G$,
consider an open subset $V$ of $X$ and obtain ${H}^U$ in $\C' (X,Y)$,
where $U \supset V$ is an open subset with $gG(V \cap M) \subset H(U \cap M)$.
Let $V'' \supset V$ be another open subset of $X$ with
$gG (V'' \cap M) \supset H(U'' \cap M) \cap gG(M)$.
Now we have the induced epimorphism
\[
g \colon G(V'' \cap M) \cap g^{-1} {H}^U \longrightarrow
{H}^U.
\]
This extends to another epimorphism
\[
g' \colon f {F}^{V''} + G(V'' \cap M) \cap g^{-1} {H}^U
\longrightarrow {H}^U
\]
with $\ker (g') = {F}^{U''}$ in $\C' (X,Y)$.
We define
\[
\mathcal{G}^V = f {F}^{U''} + G(V'' \cap M) \cap g^{-1} {H}^U.
\]
From parts (1) and (3) of Theorem \refT{lninpres}, ${G}^V$ is in $\C' (X,Y)$.
The choice for ${G}_V$ is made similarly.
\end{proof}

\begin{CorRef}{KLSDA2}
If the compactification $(X,Y)$ is small at infinity, there is a commutative square of exact inclusions
\[
\begin{CD}
\mathcal{B} (M) @>{\iota}>> \B (M) \\
@V{\kappa}VV @VV{\overline{\kappa}}V \\
\mathcal{C} (X,Y) @>{\overline{\iota}}>> \C (X,Y)
\end{CD}
\]
\end{CorRef}

\begin{proof}
The functor $\overline{\kappa}$ is the composition of exact inclusions $\B (M) \to \B'' (M)$ and the evident $\B'' (M) \to \C (X,Y)$.
\end{proof}

\SecRef{The Localization Homotopy Fibration}{LocNM}

In order for $\C (X,Y)$ to have good localization and excision properties, it is necessary to impose some geometric conditions on the pair $(X,Y)$.  In the rest of the paper we restrict to the following class of compactifications.

\begin{DefRefName}{GOOD}{Good compactification}
A compactification pair $(X,Y)$ is \textit{good} if
\begin{enumerate}
\item $X$ is metrizable, and
\item the metric on $M = X \setminus Y$ is small at infinity in the sense of Definition \refD{SmallInf}.
\end{enumerate}
\end{DefRefName}

We have already used the notion of small at infinity to relate the boundedly controlled theory $\B (M)$ from \cite{gCbG:02} to $\C (X,Y)$ in Corollary \refC{KLSDA2}.
This is the same class that was considered by Carlsson--Pedersen in \cite{gCeP:93}.

The following are some consequences of $(X,Y)$ being good.

\begin{PropRef}{UYQWEXR}
Given a closed subset $T$ of an open subset $U$ of $Y$, there is a closed subset $C$ of $X$ such that $T$ is in the interior of $C$ and $C \cap Y \subset U$.
\end{PropRef}

\begin{proof}
Since $X$ is metrizable and compact, there is a number $d$ such that $T_X [d] \cap Y$ is contained in $U$.
Here the enlargement $T_X [d]$ is taken with respect to a metric in $X$.  Clearly, for small enough $d$ the subset $T_X (d) = \{ x \in X \mid \dist (x, T) < d \}$ is an open neighborhood of $T$.
\end{proof}

\begin{PropRef}{UYQ44}
Given a closed subset $C$ of $X$, there is a subset $P$ of $X$ such that $C \subset P$, $C \cap Y = \overline{P} \cap Y$,
and $P \cap M$ is an open neighborhood of $C \cap M$.
\end{PropRef}

\begin{proof}
For any number $d > 0$, $P$ can be taken to be $C_M (d) \cup (C \cap Y)$.
Here $C_M (d)$ is the open enlargement with respect to the metric in $M$.
The equality $\overline{C_M(d)} \cap Y = C \cap Y$ holds because $C_M [d] \cap Y = C_M (d) \cap Y = C \cap Y$ by
the fact that the metric on $M$ is small at infinity.
\end{proof}

Let $F$ be an object of $\C (X,Y)$ and $Z$ be a subset of $Y$.

\begin{DefRef}{Germs}
We define $\C (X,Y)_{<Z}$ as the full subcategory of $\C (X,Y)$ on
objects supported near $Z$.
\end{DefRef}

\begin{PropRef}{Serre}
$\C (X,Y)_{<Z}$ is a Grothendieck subcategory of $\C (X,Y)$.
\end{PropRef}

\begin{proof}
To check closure under exact extensions, let
\[
F' \xrightarrow{\ f\ } F \xrightarrow{\ g\ } F''
\]
be an exact sequence in $\C (X,Y)$, where $F'$ and $F''$ are supported near $Z \subset Y$.
Therefore, there are closed subsets $C'$, $C'' \subset Y$ such that $C' \cap Y \subset Z$, $C'' \cap Y \subset Z$, and
$F'(M) \subset F'(C' \cap M)$, $F''(M) \subset F''(C'' \cap M)$.
Now $F (M) = I + N$ where $I = \im (f)$ and $N$ is any submodule of $F$ with $g(N) = F'' (M)$.
Since $g$ is continuously bicontrolled,
if $U = X \setminus C''$ then there is an open subset $V$ of $X$ such that $U \cap Y = V \cap Y$ and
$G(M \setminus U) = G(C'' \cap M) \subset fF(M \setminus V)$.
Let $D'' = X \setminus V$.
Similarly, since $f$ is continuously controlled, there is a closed subset $D'$ of $X$ with $D' \cap Y = C' \cap Y$ and
\[
fF' (M) = fF' (C' \cap M) \subset F(D' \cap M).
\]
The set $C = D' \cup D''$ is a closed subset of $X$ with the property $C \cap Y \subset Z$, and $F(C \cap M)$ contains $I + N = F$.
This shows that $F$ is supported near $Z$.

Given an admissible monomorphism $f \colon F' \to F$ in $\C (X,Y)$, suppose $F$ is supported near $Z$.
So there is a closed subset $C \subset X$ such that $C \cap Y \subset Z$ and $F(M) \subset F(C \cap M)$.
Notice that $fF' (M) \subset F(M) \subset F(C \cap M)$.
Since $f$ is continuously bicontrolled, we have a subset $C'$ with $C' \cap Y \subset Z$ and $F' (M) \subset F' (C' \cap M)$.

Closure of $\C (X,Y)_{<Z}$ under admissible quotients is similar.
\end{proof}

Let $\Z$ be the subcategory $\C (X, Y)_{<Z}$ of
$\C (X, Y)$. Let the class
of \textit{weak equivalences} $\Sigma$ in $\C (X, Y)$ consist of all finite
compositions of admissible monomorphisms with cokernels in $\Z$ and
admissible epimorphisms with kernels in $\Z$.  We want to show that the
class $\Sigma$ admits a calculus of right fractions.

\begin{LemRef}{JHQASF}
Suppose $G$ is a graded module in $\C (X,Y)$ with a grading given by $\mathcal{G}^{\bullet}$ and $\mathcal{G}_{\bullet}$.  Let $F$ be a submodule which is split at infinity with respect to the standard filtration.  Then the assignments ${F}^V = F \cap {G}^V$ and ${F}_V = F \cap {G}_V$ give a grading of $F$.
\end{LemRef}

\begin{proof}
Of course, $F (V \cap M) = F \cap G(V \cap M) \subset F \cap {G}^{V} = {F}^{V}$.
On the other hand, $G^V \subset G(V' \cap M)$ for some $V'$ as in property (1) of Definition \refD{graded2}, so ${F}^V \subset F \cap G(V' \cap M) = F(V' \cap M)$.
Property (2) is checked similarly.

Consider the inclusion $i \colon F \to G$, and take the quotient $q \colon G \to H$.
Both $F$ and $G$ are split and insular at infinity, so $H$ is split and insular by parts (4) and (5) of Theorem \refT{lninpres}, with respect to the quotient filtration.
We define ${H}^{\bullet}$ and ${H}_{\bullet}$ as the images $q{G}^{\bullet}$ and $q{G}_{\bullet}$ and give $H^V$ and $H_V$ the standard filtration in $H$.  Then both of them are split at infinity as the images of $G^V$ and $G_V$, which are split at infinity, and insular since $H$ is insular.  Now the kernel of $q \vert \colon {G}^V \to {H}^V$,
which is $F \cap {G}^V$, is split by part (5) of Theorem \refT{lninpres}.
Since $F$ is insular at infinity, ${F}^V$ is also insular at infinity.
Similarly, $F_V$ is split and insular at infinity.
So the assignments ${F}^{\bullet}$ and ${F}_{\bullet}$ give a grading for $F$.
\end{proof}

This result can be promoted to the following statement.

\begin{LemRef}{HUVZFMN}
Suppose $F$ is the kernel of a continuously bicontrolled epimorphism $g \colon G \to H$ in $\C' (X,Y)$.  If $G$ is a graded module in $\C (X,Y)$ and $F$ is split at infinity then both $H$ and $F$ are in $\C (X,Y)$.
\end{LemRef}

\begin{proof}
The grading for $H$ is given by ${H}^V = g {G}^U$ and ${H}_V = g {G}_{W}$ where $U$ and $W$ are determined by the assumption on $g$, as in Definition \refD{ContCat2}.
Each ${H}^V$ and ${H}_V$ is split and insular at infinity by the arguments in the proof of Lemma \refL{JHQASF}.
The inclusions $H(V \cap M) \subset gG(U \cap M) \subset g {G}^U = {H}^V$,
$H^V = g \mathcal{G}^U \subset g {G}(V' \cap M) \subset H(U' \cap M)$, and the similar inclusions verifying property (2) from Definition \refD{graded2}, show that ${H}^{\bullet}$ and ${H}_{\bullet}$ give a grading for $H$.
The same argument shows that the assignments ${F}^V = F \cap {G}^U$ and ${F}_V = F \cap {G}_{W}$ gives a grading for $F$.
\end{proof}

\begin{DefRef}{rfil}
A Grothendieck subcategory $\Z \subset \C$ is \textit{right
filtering} if each morphism $f \colon F \to G$ in $\C$, where $G$
is an object of $\Z$, factors through an admissible epimorphism $e
\colon F \to \overline{G}$, where $\overline{G}$ is in $\Z$.
\end{DefRef}

\begin{LemRef}{SubSubLem}
The subcategory $\Z = \C (X, Y)_{<Z}$ of $\C = \C (X, Y)$ is right filtering.
\end{LemRef}

\begin{proof}
Since $G$ is an object of $\Z$, there is a closed subset $C$ such that $C \cap Y \subset Z$ and $G(M) \subset G(C\cap M)$.
Using that $G$ is insular at infinity, there exists an open subset $U \subset X \setminus C$ such that $G(U \cap M) = 0$.
Then there is an open subset $V \subset U$ such that $V \cap Y = Y \setminus C$ and $fF(V \cap M) \subset G(U \cap M)$,
so $fF(V \cap M) =0$.
Choose $V'' \subset V$ with $V'' \cap Y = V \cap Y$ and $F_V \subset F(V \cap M)$ guaranteed by the grading.
By Lemma \refL{HUVZFMN}, $F_V$ is in $\C (X,Y)$ and the cokernel $\overline{G}$ of the inclusion $F_V \to F$
is also an object of $\C (X,Y)$.
In fact, $\overline{G}$ is supported near $C \cap Y \subset Z$, so $\overline{G}$ is in $\Z$.
Since $fF(V \cap M) =0$,
the required factorization is the right square in the map between the two exact sequences
\[
\begin{CD}
F_V @>>> F @>{j'}>> \overline{G}\\
@V{i}VV @VV{=}V @VVV \\
K @>{k}>> F @>{f}>> G
\end{CD}
\]
\end{proof}

\begin{CorRef}{WeakEqs}
The class $\Sigma$ admits a calculus
of right fractions.
\end{CorRef}

\begin{proof}
This follows from Lemma \refL{SubSubLem}, see Lemma 1.13 of \cite{mS:03}.
\end{proof}

\begin{DefRef}{Quot}
The
\define{quotient category}
$\mathbf{C}/\mathbf{Z}$ is the localization $\C [\Sigma^{-1}]$.
\end{DefRef}

It is clear that the quotient $\mathbf{C}/\mathbf{Z}$ is an
additive category, and $P_{\Sigma}$ is an additive functor.
In fact, we have the following.

\begin{ThmRef}{ExLocStr}
The short sequences in $\mathbf{C}/\mathbf{Z}$ which are
isomorphic to images of exact sequences from $\C$ form a Quillen
exact structure.
\end{ThmRef}

\begin{proof}
This will follow from Proposition 1.16 of Schlichting \cite{mS:03}.
Since $\Z$ is right filtering by Lemma \refL{SubSubLem}, it remains
to check that the subcategory $\Z$ is \textit{right s-filtering} in $\C$,
that is to show that if $f \colon F \to G$ is an admissible
monomorphism with $F$ in $\Z$ then there exist $E$ in $\Z$ and an
admissible epimorphism $e \colon G \to E$ such that the composition $ef$ is an
admissible monomorphism.

Assume that $F$ is supported near $Z$, then the image $f(F)$ is also supported near $Z$ with the appropriate choice of a closed subset $C$ such that $C \cap Y \subset Z$.
If $V$ is an open subset with $V \cap Y = Y \backslash C$, choose $G_V$ given by the grading of $G$ so that $G_V$ is in $\C$ and $G(C) \cap G_V = 0$.
Let $e \colon G \to E$ be the cokernel of the inclusion $G_V \to G$.
Since $f (F) \cap G_V = 0$, $ef$ is an admissible
monomorphism.
\end{proof}

\begin{NotRef}{GERMS}
In the case when $\Z = \C (X, Y)_{<Z}$ and $\C = \C (X, Y)$, we will use the notation $\C (X, Y)^{>Z}$ for the exact quotient $\mathbf{C}/\mathbf{Z}$.
\end{NotRef}

The Quillen algebraic $K$-theory is a functor from the category of exact categories and exact functors to connective spectra.

Now we have available the following localization theorem.

\begin{ThmRefName}{Schlichting}{Theorem 2.1 of Schlichting \cite{mS:03}}
Let $\A$ be an idempotent complete right s-filtering subcategory
of an exact category $\E$.  Then the sequence of exact categories
$\A \to \E \to {\E}/{\A}$ induces a homotopy fibration of $K$-theory
spectra
\[
K(\A) \longrightarrow K(\E) \longrightarrow K({\E}/{\A}).
\]
\end{ThmRefName}

Applying Schlichting's theorem to the subcategory $\Z = \C (X, Y)_{<Z}$ of $\C = \C (X, Y)$, and denoting by $G (X, Y)$ the $K$-theory of $\C (X, Y)$, we obtain the basic localization theorem.

\begin{ThmRefName}{Sch}{Localization}
There is a homotopy fibration
\[
G (X, Y)_{<Z} \longrightarrow G (X, Y) \longrightarrow G (X, Y)^{>Z}.
\]
\end{ThmRefName}

This is easily generalized to the following localization result.  We leave the details to the reader.

\begin{CorRef}{Sch2}
Suppose $W$ is an open subset of $Y$, and $Z$ is a subset of $W$.
Then there is a homotopy fibration
\[
G (X, Y)_{<Z} \longrightarrow G (X, Y)_{<W} \longrightarrow G (X, Y)^{>Z}_{<W}.
\]
\end{CorRef}

\SecRef{Controlled Excision Theorems}{LocK}

The controlled excision theorem will be stated in terms of nonconnective $K$-theory spectra which are specific nonconnective deloopings of the $K$-theory of our exact categories.
The construction of the deloopings and the proof of controlled excision require the context of
Waldhausen categories with cofibrations and weak equivalences and their $K$-theory.
We will assume that the reader is familiar with the terminology and some results reviewed in
section 4 of \cite{gCbG:02}, which are all standard.
We only briefly recall the statement of the main result used in this paper.

Let $\D$ be a small Waldhausen category with respect to two
categories of weak equivalences $\bsv(\D) \subset \w(\D)$ with a
cylinder functor $T$ both for $\vD$ and for $\wD$ satisfying the
cylinder axiom for $\wD$. Suppose also that $\w(\D)$ satisfies the
extension and saturation axioms.

We define $\vDw$ to be the full
subcategory of $\vD$ whose objects are $F$ such that $0 \to F \in
\w(\D)$. Then $\vDw$ is a small Waldhausen category with
cofibrations $\co (\Dw) = \co (\D) \cap \Dw$ and weak
equivalences $\bsv (\Dw) = \bsv (\D) \cap \Dw$. The cylinder
functor $T$ for $\vD$ induces a cylinder functor for $\vDw$.  Since
$T$ satisfies the cylinder axiom, the induced functor does so
too.

\begin{ThmRefName}{ApprThm}{Approximation theorem}
Let $E \colon \D_1 \rightarrow \D_2$ be an exact functor between
two small saturated Waldhausen categories. It induces a map of
$K$-theory spectra
\[
  K (E) \colon K (\D_1) \longrightarrow K (\D_2).
\]
Assume that $\D_1$ has a cylinder functor satisfying the cylinder
axiom. If $E$ satisfies two conditions:
\begin{enumerate}
\item a morphism $f \in \D_1$ is in $\w(\D_1)$
if and only if $E (f) \in \D_2$ is in $\w(\D_2)$,
\item for any object $D_1 \in \D_1$ and any morphism
$g \colon E(D_1) \to D_2$ in $\D_2$, there is an object $D'_1 \in
\D_1$, a morphism $f \colon D_1 \to D'_1$ in $\D_1$, and a weak
equivalence $g' \colon E(D'_1) \rightarrow D_2 \in \w(\D_2)$ such
that $g = g' E(f)$,
\end{enumerate}
then $K (E)$ is a homotopy equivalence.
\end{ThmRefName}

In the context of an arbitrary additive category, a sequence of morphisms
\[
E^{\subdot} \colon \quad 0 \longrightarrow E^1 \xrightarrow{\ d_1
\ } E^2 \xrightarrow{\ d_2 \ }\ \dots\ \xrightarrow{\ d_{n-1} \ }
E^n \longrightarrow 0
\]
is called a \textit{(bounded) chain complex} if the compositions
$d_{i+1} d_i$ are the zero maps for all $i = 1$,\dots, $n-1$.  A
\textit{chain map} $f \colon F^{\subdot} \to E^{\subdot}$ is a
collection of morphisms $f^i \colon F^i \to E^i$ such that $f^i
d_i = d_i f^i$.  A chain map $f$ is \textit{null-homotopic} if
there are morphisms $s_i \colon F^{i+1} \to E^i$ such that $f = ds
+ sd$.  Two chain maps $f$, $g \colon F^{\subdot} \to E^{\subdot}$
are \textit{chain homotopic} if $f-g$ is null-homotopic. Now $f$
is a \textit{chain homotopy equivalence} if there is a chain map
$h \colon E^i \to F^i$ such that the compositions $fh$ and $hf$
are chain homotopic to the respective
identity maps.

The Waldhausen structures on categories of bounded chain complexes
are based on homotopy equivalence as a weakening of the notion of
isomorphism of chain complexes.

A sequence of maps in an exact category is called \textit{acyclic}
if it is assembled out of short exact sequences in the sense that
each map factors as the composition of the cokernel of the
preceding map and the kernel of the succeeding map.

It is known that the class of acyclic complexes in an exact
category is closed under isomorphisms in the homotopy category if
and only if the category is idempotent complete, which is also
equivalent to the property that each contractible chain complex is
acyclic, cf.\ \cite[sec.\ 11]{bK:96}.

Given an exact category $\E$, there is a standard choice for the
Waldhausen structure on the category $\E'$ of bounded
chain complexes in $\E$ where the degree-wise admissible
monomorphisms are the cofibrations and the chain maps whose
mapping cones are homotopy equivalent to acyclic complexes are the
weak equivalences $\boldsymbol{v}(\E')$.
We denote this Waldhausen structure by $\vE'$.

The following result is verified in \cite{gCbG:02}.

\begin{PropRef}{FibThApp2}
The category $\vE'$ is a Waldhausen category satisfying the
extension and saturation axioms and has cylinder functor
satisfying the cylinder axiom.
\end{PropRef}

This associates to the exact continuously controlled category $\C (X,Y)$ the Waldhausen category $\vC'$.

Another choice for the Waldhausen structure on the category of bounded chain complexes
in $\C$ can be associated to a given subset $Z$ of the boundary $Y$.
The new weak equivalences $\w (\C')$ are the chain
maps whose mapping cones are homotopy equivalent to acyclic
complexes in the quotient $\mathbf{C}/\mathbf{Z}$.

\begin{CorRef}{FibSt}
The categories $\vC'$ and $\wC'$ are Waldhausen categories
satisfying the extension and saturation axioms and have cylinder
functors satisfying the cylinder axiom.
\end{CorRef}

The $K$-theory functor from the category of small Waldhausen
categories $\D$ and exact functors to connective spectra is
defined in terms of $S_{\subdot}$-construction as in Waldhausen
\cite{fW:85}. It extends to simplicial categories $\D$ with
cofibrations and weak equivalences and inductively delivers the
connective spectrum $n \mapsto \vert \bfw S_{\subdot}^{(n)} \D
\vert$. We obtain the functor assigning to $\D$ the connective
$\Omega$-spectrum
\[
K (\D) = \Omega^{\infty} \vert \bfw S_{\subdot}^{(\infty)} \D
\vert = \colim{n \ge 1} \Omega^n \vert \bfw S_{\subdot}^{(n)} \D
\vert
\]
representing the Waldhausen algebraic $K$-theory of $\D$. For
example, if $\D$ is the additive category of free finitely
generated $R$-modules with the canonical Waldhausen structure,
then the stable homotopy groups of $K (\D)$ are the usual
$K$-groups of the ring $R$.  In fact, there is a general
identification of the two theories. Recall that for an exact
category $\E$, the category $\E'$ of bounded chain complexes has the Waldhausen
structure $\vE'$.
We know the following from \cite{gCbG:02}.

\begin{ThmRef}{Same}
The Quillen $K$-theory of an exact category $\E$ is equivalent to
the Waldhausen $K$-theory of $\vE'$.
\end{ThmRef}

We now apply a construction from \cite{gCbG:02} to the exact category $\C = \C (X,Y)$.
It requires an embedding of $\C$ in a pseudoabelian exact category, for that we choose the inclusion of $\C$ in $\U (X,Y)$.
The suspension $S \mathbf{C}$ of $\C$ is defined as a certain exact quotient in \cite{gCbG:02}, and there is a map $K (\C) \to \Omega K (S \mathbf{C})$ which is a
weak equivalence in positive dimensions. Iterations of this
construction give weak equivalences $K (S^k \mathbf{C}) \to \Omega
K (S^{k+1} \mathbf{C})$.

\begin{DefRef}{NonconnKw}
The \textit{nonconnective continuously controlled $G$-theory}
is the spectrum
\[
\Gnc (X,Y) \overset{ \text{def} }{=} \hocolim{k > 0}
\Omega^{k} K (S^k \mathbf{C}).
\]
\end{DefRef}

Since $\Gnc (X,Y)$ is an $\Omega$-spectrum in positive
dimensions, the positive homotopy groups of $\Gnc (X,Y)$ coincide
with those of $G (X,Y)$ as desired.  In fact, the class groups
are also isomorphic as $\C$ is idempotent complete.

If $Z$ is a subset of $Y$, recall the
fibration
\[
G (X,Y)_{<Z} \longrightarrow G (X,Y) \longrightarrow G (X,Y)^{>Z}.
\]
The terminal spectrum is the $K$-theory of the quotient ${\C}/{\Z}$.
Notice that there is a map
\[
K ({\C}/{\Z}) \longrightarrow \Omega K ({S\C}/{S\Z})
\]
which is a weak equivalence in positive dimensions by the Five
Lemma. If one defines
\[
\Knc ({\C}/{\Z}) = \hocolim{k}
\Omega^{k} K ({S^k \C}/{S^k \Z}),
\]
there is an induced fibration
\[
\Gnc (X,Y)_{<Z} \longrightarrow \Gnc (X,Y) \longrightarrow \Gnc (X,Y)^{>Z}.
\]

The following theorem is the main result of the paper.

\begin{ThmRefName}{ExCC}{Continuously Controlled Excision}
Let $(X,Y)$ be a good compactification of a metric space $M = X \setminus Y$ by attaching the boundary $Y$.
For any pair of open subsets $U_1$, $U_2 \subset Y$, there is a
homotopy pushout diagram
\[
\begin{CD}
\Gnc (X,Y)_{<U_1 \cap U_2} @>>> \Gnc (X,Y)_{<U_1} \\
@VVV @VVV \\
\Gnc (X,Y)_{<U_2} @>>> \Gnc (X,Y)_{<U_1 \cup U_2}
\end{CD}
\]
\end{ThmRefName}

The proof of Theorem \refT{ExCC} will require a sequence of lemmas.

\begin{LemRef}{CharWE}
If $f^{\subdot} \colon F^{\subdot} \to G^{\subdot}$ is a
degreewise admissible monomorphism with cokernels in $\Z$ then
$f^{\subdot}$ is a weak equivalence in $\boldsymbol{v}({\C}/{\Z})'$.
\end{LemRef}

\begin{proof}
We need to see that the mapping cone $Cf^{\subdot}$ is the zero
complex in the homotopy category of ${\C}/{\Z}$.
But $Cf^{\subdot}$ is weakly equivalent to the cokernel of
$f^{\subdot}$, by Lemma 11.6 of \cite{bK:96}, which is zero in
${\C}/{\Z}$.
\end{proof}

\begin{LemRef}{CharWE2}
Suppose $U_1$ and $U_2$ form an open cover of $Y$.
Given a continuously controlled homomorphism $g \colon F \to G$, suppose $F$ is an object of $\C (X,Y)_{<U_1}$.
Then there is an admissible subobject $\iota \colon G' \to G$ with $G'$ such that $g$ factors through $\iota$, $G'$ is an object of $\C (X,Y)_{<U_1}$, and the quotient $G/G'$ is an object of $\C (X,Y)_{<U_2}$.
\end{LemRef}

\begin{proof}
Suppose $F$ is supported on a closed subset $C'$ of $X$ with $C' \cap Y \subset Z$.
Let $C''$ be a closed neighborhood of $U_1 \backslash U_2$ such that $C'' \cap Y \subset U_1$ from Proposition \refP{UYQWEXR}.
Let $C = C' \cup C''$.
The complement $X \setminus C$ is open, so one can use the fact that $g$ is continuously controlled at infinity to find an open subset $W$ such that $W \cap C = \emptyset$, $gF(W \cap M) \subset G(M \setminus C)$, and $g F(C \cap M) \subset G(X \setminus W \cap M)$.
Choose $V = X \setminus W$ and let $G' = G^V$ be the submodule of $G$ given by the grading.
It is an admissible subobject by Lemma \refL{HUVZFMN}.
\end{proof}

The following is a slight relativization of the last result.

\begin{LemRef}{CharWE3}
Suppose $U_1$ and $U_2$ are open subsets of $Y$.
Given a continuously controlled homomorphism $g \colon F \to G$ in $\C (X,Y)_{<U_1 \cup U_2}$, suppose $F$ is an object of $\C (X,Y)_{<U_1}$.
Then there is an admissible subobject $\iota \colon G' \to G$ with $G'$ such that $g$ factors through $\iota$, $G'$ is an object of $\C (X,Y)_{<U_1}$, and the quotient $G/G'$ is an object of $\C (X,Y)_{<U_2}$.
\end{LemRef}

\begin{proof}
This proceeds along the same lines but using a variant of Proposition \refP{UYQWEXR} where $T$ is closed in $U$ but not necessarily closed in $Y$.
\end{proof}

The last two lemmas can be generalized to the following.

\begin{LemRef}{CharWE33}
Given continuously controlled homomorphisms $g \colon F \to G$ and $\psi \colon F' \to G$ in $\C (X,Y)_{<U_1 \cup U_2}$, suppose both $F$ and $F'$ are objects of $\C (X,Y)_{<U_1}$.  Then there is an admissible subobject $\iota \colon G' \to G$ such that both $g$ and $\psi$ factor through $\iota$, $G'$ is an object of $\C (X,Y)_{<U_1}$, and the quotient $G/G'$ is an object of $\C (X,Y)_{<U_2}$.
\end{LemRef}

\begin{proof}
One starts with the two constructions from the proof of Lemma \refL{CharWE2} performed independently.  Once the subsets that play the role of $W$ have been constructed for both maps, let $W$ instead denote the intersection of those open subsets.
Then simply proceed to obtain $V$ and the subobject $G' = G^V$ as before.
\end{proof}

\begin{proof}[Proof of Theorem \refT{ExCC}]
We are using the notation $\C_{1,2} = \C (X,Y)_{<U_1 \cup U_2}$, $\C_i= \C
(X,Y)_{<U_i}$ for $i=1$ or $2$, and $\C_{12} = \C (X,Y)_{<U_1 \cap U_2}$.
There is a commutative diagram
\[
\begin{CD}
K (\C_{12}) @>>> K (\C_1) @>>> K ({\C_1}/{\C_{12}}) \\
@VVV @VVV @VV{K(I)}V \\
K (\C_2) @>>> K (\C_{1,2}) @>>> K (\C_{1,2} \! /{\C_2})
\end{CD}
\]
where the rows are homotopy fibrations by Corollary \refC{Sch2}.

Here $I \colon {\C_1}/{\C_{12}} \to \C_{1,2} \! /{\C_2}$ is not an equivalence of categories as it usually is in similar
computations where Karoubi filtrations are available, cf.~\cite{mCeP:97}.
We only claim that one can apply the Approximation Theorem \refT{ApprThm} to $I \colon \wC'_1 \to \wC'_{1,2}$ and conclude that
$K (I)$ is a weak equivalence.

The first condition is easy.
To check the second condition, consider
\[
F^{\subdot} \colon \quad 0 \longrightarrow F^1 \xrightarrow{\
\phi_1\ } F^2 \xrightarrow{\ \phi_2\ }\ \dots\
\xrightarrow{\ \phi_{n-1}\ } F^n \longrightarrow 0
\]
in $\C_1$ and a chain map $g \colon F^{\subdot} \to G^{\subdot}$ to a chain complex
\[
G^{\subdot} \colon \quad 0 \longrightarrow G^1 \xrightarrow{\
\psi_1\ } G^2 \xrightarrow{\ \psi_2\ }\ \dots\
\xrightarrow{\ \psi_{n-1}\ } G^n \longrightarrow 0
\]
in $\C$.

The fact that $F^{\subdot}$ is a complex in $\C_1$ allows us to make specific choices of subsets of $X$ to support each module $F^i$.
Suppose $F^1$ is supported on a closed subset $C_1$ of $X$ with $C_1 \cap Y \subset U_1$.
We can now proceed to use Lemma \refL{CharWE2} and construct an admissible subobject $F'^1$ of $G^1$ which is in $\C_1$ and receives the restriction $g_1 \colon F^1 \to F'^1$.
We can now apply Lemma \refL{CharWE33} to the two maps $g_2 \colon F^2 \to G^2$ and $\psi_1 \colon F'^1 \to G^2$ to construct the admissible subobject $F'^2$ of $G^2$.
We define the complex $F^{\prime{i}}$ by inductive application of Lemma \refL{CharWE3} and define $\xi_{i} \colon F^{\prime{i}} \to F^{\prime{i+1}}$ to be the restrictions of
$\psi_{i}$ to $F^{\prime{i}}$.
This gives a chain subcomplex
$(F^{\prime{i}},\xi_{i})$ of $(G^{i},\psi_{i})$ in $\C_1$ with the inclusion $s \colon F^{\prime{i}} \to G^{i}$.
Notice that the choices give the induced chain map $\overline{g} \colon F^{\subdot} \to F^{\prime\subdot}$ in $\C_1$ so that $g = s \circ I(\overline{g})$.
Now $s$ is a weak equivalence of chain complexes in
$\C_{1,2} \! /{\C_2}$ by Lemma \refL{CharWE} since $C^{\subdot} = \coker (s)$ is in $\C_2$.
This implies that after applying nonconnective deloopings, the rightmost vertical arrow in the commutative diagram
\[
\begin{CD}
\Gnc (X,Y)_{<U_1 \cap U_2} @>>> \Gnc (X,Y)_{<U_1}
@>>> \Gnc (X,Y)_{<U_1}^{>U_1 \cap U_2}\\
@VVV @VVV @VVV\\
\Gnc (X,Y)_{<U_2} @>>> \Gnc (X,Y)_{<U_1 \cup U_2}
@>>> \Gnc (X,Y)_{<U_1 \cup U_2}^{>U_2}
\end{CD}
\]
is a weak equivalence, so the homotopy colimit of the square on the left is a homotopy
pushout.
\end{proof}

The case of an open covering $\{ U_1, U_2 \}$ of $Y$ gives the following consequence which was stated in the introduction.

\begin{CorRef}{ExCCC}
Let $(X,Y)$ be a good compactification of a metric space $M = X \setminus Y$ by attaching the boundary $Y$.
There is a
homotopy pushout square
\[
\begin{CD}
\Gnc (X,Y)_{<U_1 \cap U_2} @>>> \Gnc (X,Y)_{<U_1} \\
@VVV @VVV \\
\Gnc (X,Y)_{<U_2} @>>> \Gnc (X,Y)
\end{CD}
\]
\end{CorRef}

We finally give examples of the geometric situations where the given metric space has a good compactification.

\begin{CorRef}{ExCCCC}
Suppose $M$ is a complete metric space whose metric is either hyperbolic in the sense of Gromov or is CAT(0) and let $X$ be the compactification of $M$ is by attaching the ideal boundary $Y$.
Then the commutative diagrams in Theorem \refT{ExCC} and Corollary \refC{ExCCC} are homotopy pushouts.
\end{CorRef}

\begin{proof}
The facts that the ideal compactifications of hyperbolic and CAT(0) spaces are good can be found in
\cite{eGpH:90} and \cite{mBaH:99}.
\end{proof}


\begin{thebibliography}{999}

\bibitem{dAfCsFeP:94}
{D.R.~Anderson, F.X.~Connolly, S.C.~Ferry, and E.K.~Pedersen},
{\it Algebraic $K$-theory with continuous control at infinity},
J. Pure Appl. Alg. \textbf{94} (1994), 25--47.

\bibitem{aB:03}
{A.~Bartels}, \textit{Sqeezing and higher algebraic \textit{K}-theory}, \textit{K}-theory {\bf 28},
(2003), 19--37.

\bibitem{BFJR:04}
{A.~Bartels, T.~Farrell, L.~Jones, H.~Reich},
{\it On the Isomorphism Conjecture in algebraic $K$-theory},
Topology,  \textbf{43} (2004), 157--213.

\bibitem{mBaH:99}
M. Bridson and A. Haefliger, \textit{Metric spaces of non-positive curvature}, Springer-Verlag, 1999.

\bibitem{mCeP:97}
{M.~Cardenas and E.K.~Pedersen}, {\it On the Karoubi filtration of
a category}, $K$-theory, {\bf 12} (1997), 165--191.

\bibitem{gC:95}
{G. Carlsson}, {\it Bounded $K$-theory and the assembly map in
algebraic $K$-theory}, in {\it Novikov conjectures, index theory
and rigidity}, {\it Vol. 2} (S.C.~Ferry, A.~Ranicki, and
J.~Rosenberg, eds.), Cambridge U. Press (1995), 5--127.

\bibitem{gCbG:04}
{G. Carlsson and B. Goldfarb}, \textit{The integral K-theoretic Novikov conjecture for
groups with finite asymptotic dimension}, Inventiones Math. {\bf
157} (2004), 405--418.

\bibitem{gCbG:03}
\bysame, \textit{On homological coherence of
discrete groups}, J.~Algebra {\bf 276} (2004), 502--514.

\bibitem{gCbG:02}
\bysame, \textit{Controlled algebraic \textit{G}-theory, I},
J. Homotopy Relat. Struct. \textbf{6} (2011), 119–-159.

\bibitem{gCbG:13}
\bysame, \textit{Algebraic $K$-theory of geometric groups}, preprint, 2013.
\texttt{arXiv:1305.3349}

\bibitem{gCeP:93}
{G. Carlsson and E.K. Pedersen}, {\it Controlled algebra and the
Novikov conjecture for $K$- and $L$-theory}, Topology {\bf 34}
(1993), 731--758.

\bibitem{gCeP:97}
\bysame, {\it \v{C}ech homology and the Novikov conjectures for
$K$- and $L$-theory}, Math. Scand. \textbf{82} (1998), 5--47.

\bibitem{eGpH:90}
E. Ghys and P. de la Harpe (ed.), \textit{Sur les Groupes Hyperboliques d'apr`es Mikhael Gromov},
Progr. Math. \textbf{83}, Birkh\"{a}user, 1990.

\bibitem{bG:97}
{B. Goldfarb}, {\it Novikov conjectures for arithmetic groups
with large actions at infinity}, $K$-theory {\bf 11} (1997),
319--372.

\bibitem{bG:98}
{\bysame}, {\it Novikov conjectures and relative hyperbolicity},
Math. Scand. \textbf{85} (1999), 169--183.

\bibitem{bG:99}
{\bysame}, \textit{Large scale topology of arithmetic groups and the strong Novikov conjecture},
Topology Appl. (2003).

\bibitem{bK:96}
{B.~Keller}, \textit{Derived categories and their uses}, in
Handbook of Algebra, Vol.~1 (M.~Hazewinkel, ed.), 1996, Elsevier
Science, 671--701.

\bibitem{tL:07}
{T.K.~Lance}, \textit{Continuous control and Novikov conjectures in exact non-split categories},
Ph.D. thesis, SUNY--Albany, 2007.

\bibitem{eP:95}
{E.K. Pedersen},
{\it Bounded and continuous control},
in {\it Novikov conjectures, index theory and rigidity},
{\it Vol. 2}
(S.C.~Ferry, A.~Ranicki, and J.~Rosenberg, eds.),
Cambridge U. Press (1995), 277--284.

\bibitem{ePcW:85}
{E.K. Pedersen and C. Weibel}, {\it A nonconnective delooping of
algebraic $K$-theory}, in {\it Algebraic and geometric topology}
(A.~Ranicki, N.~Levitt, and F.~Quinn, eds.), Lecture Notes in
Mathematics {\bf 1126}, Springer-Verlag (1985), 166--181.

\bibitem{ePcW:89}
{\bysame}, {\it $K$-theory homology of spaces}, in {\it Algebraic
topology} (G.~Carlsson, R.L.~Cohen, H.R.~Miller, and
D.C.~Ravenel, eds.), Lecture Notes in Mathematics {\bf 1370},
Springer-Verlag (1989), 346--361.

\bibitem{dRrTgY:11}
{D.~Ramras, R.~Tessera, and G.~Yu},
\textit{Finite decomposition complexity and the integral Novikov conjecture for higher algebraic K-theory},
preprint (2011). \texttt{arXiv:1111.7022}

\bibitem{dR:04}
{D.~Rosenthal},
\textit{Splitting with continuous control in algebraic K-theory},
$K$-theory {\bf 32} (2004), 139--166.

\bibitem{mS:03}
{M.~Schlichting}, \textit{Delooping the $K$-theory of exact categories},
Topology {\bf 43} (2004), 1089--1103.

\bibitem{rTtT:90}
{R.W. Thomason and Thomas Trobaugh}, {\it Higher algebraic
$K$-theory of schemes and of derived categories}, in {\it The
Grothendieck Festschrift}, {\it Vol. III}, Progress in
Mathematics {\bf 88}, Birkh\"{a}user (1990), 247--435.

\bibitem{fW:85}
{F. Waldhausen}, {\it Algebraic $K$-theory of spaces}, in {\it
Algebraic and geometric topology} (A. Ranicki, N. Levitt, and F.
Quinn, eds.), Lecture Notes in Mathematics {\bf 1126},
Springer-Verlag (1985), 318--419.

\end{thebibliography}
\end{document}